\documentclass[english,11pt]{article}
\usepackage{amsthm}
\usepackage{fullpage}
\usepackage[margin=1in]{geometry}
\usepackage{amsmath}
\usepackage{mathrsfs}
\usepackage{amssymb}
\usepackage{subfigure}
\usepackage{verbatim}
\usepackage{caption}
\usepackage{color}
\usepackage{enumerate}

\usepackage{array}
\usepackage{graphicx}
\def\qed{\hfill \ifhmode\unskip\nobreak\fi\quad\ifmmode\Box\else$\Box$\fi\\ }
\usepackage{subfigure}

\newcommand{\beq}[1]{\begin{equation}\label{#1}}
\newcommand{\eeq}{\end{equation}}
\newcommand{\blem}[1]{\begin{lemma}\label{#1}}
\newcommand{\elem}{\end{lemma}}
\newcommand{\bth}[1]{\begin{theorem}\label{#1}}
\newcommand{\enth}{\end{theorem}}
\newcommand{\brem}[1]{\begin{remark}\label{#1}}
\newcommand{\erem}{\end{remark}}

\newtheorem{theorem}{Theorem}
\newtheorem{lemma}[theorem]{Lemma}
\newtheorem{prop}[theorem]{Proposition}
\newtheorem{corollary}[theorem]{Corollary}

\newtheorem{conjecture}[theorem]{Conjecture}
\theoremstyle{definition}
\newtheorem{defn}{Definition}

\begin{document}

\title{Injective edge-coloring of  graphs with given maximum degree}
\author{
Alexandr Kostochka\thanks{Department of Mathematics, University of Illinois at Urbana-Champaign, Urbana, IL 61801, USA, and Sobolev Institute of Mathematics, Novosibirsk 630090, Russia. Email: {\tt kostochk@math.uiuc.edu.}  Partially supported by NSF grant DMS-1600592,   NSF RTG Grant DMS-1937241,
  ANR project HOSIGRA (ANR-17-CE40-0022) and  grant 19-01-00682 of the Russian Foundation for Basic Research.}
\and 
{Andr\'e Raspaud}
\thanks{LaBRI (Universit\'e de Bordeaux), 351 cours de la Lib\'eration, 33405 Talence Cedex, France. {E-mail: andre.raspaud@labri.fr. This work is supported by the ANR project HOSIGRA (ANR-17-CE40-0022).}}
\and
Jingwei Xu\thanks{Department of Mathematics, University of Illinois at Urbana-Champaign, Urbana, IL 61801, USA. Email:{\tt jx6@illinois.edu}.
Partially supported   by Arnold O. Beckman Campus Research Board Award RB20003 of the University of Illinois at Urbana-Champaign.}}


\maketitle
\begin{abstract}
A coloring of edges of a graph $G$ is {\em injective} if for any two distinct edges $e_1$ and $e_2$, the colors of 
$e_1$ and $e_2$ are distinct if they are at distance $1$ in $G$ or in a common triangle. 
 Naturally, the
 injective chromatic index of $G$,  $\chi'_{\rm {inj}}(G)$, is the minimum number of colors needed for an injective edge-coloring of $G$. 
 We study how large can be the
 injective chromatic index of $G$ in terms of maximum degree of $G$ when we have restrictions on girth and/or chromatic number of~$G$. We also compare our bounds with analogous bounds on the strong chromatic index.
\vspace{0.5cm}

\noindent
 {\small{\em Mathematics Subject Classification}: 05C15, 05C35.}\\
 {\small{\em Key words and phrases}:  Injective coloring,  injective edge coloring, maximum degree.}
\end{abstract}

\section{Introduction}
\subsection{Notation and definitions}

For a positive integer $n$, we denote by $[n]$ the set $\{1,\ldots,n\}$.
By $\Delta(G)$ we denote the maximum degree of a graph $G$ and by $\alpha(G)$ -- the {\em independence number} of $G$.

For disjoint subsets $A$ and $B$ of vertices in a graph $G$, let $E_G(A,B)$ denote the set of edges in $G$ with one end in $A$ and one in $B$.
Also $G[A]$ denotes the subgraph of $G$ induced by $A$.
By $N_G(v)$ we denote the neighborhood of vertex $v$ in graph $G$, and let $d_G(v)=|N_G(v)|$. When the graph $G$ is clear from the context, we can drop the subscript. 

A vertex coloring of a graph $G$ is {\em injective} if for every vertex $v$ of $G$, all the neighbors of $v$ have different colors. In other words, in 
injective coloring two vertices $u$ and $v$ must have distinct colors if there is a $u,v$-path of length exactly $2$.
The {\em injective chromatic number}, $\chi_{\rm {inj}}(G)$, of a graph $G$ is the minimum $k$ such that $G$ admits injective coloring with $k$ colors.

Similarly, an edge coloring of a graph $G$ is {\em injective} if any
 two edges $e$ and $f$  that are at distance exactly $1$ in $G$ or are in a common triangle have distinct colors.
The {\em injective chromatic index} of $G$,  $\chi'_{\rm {inj}}(G)$, is the minimum number of colors needed for an injective edge coloring of $G$. 

Note that an injective edge coloring is not necessarily a proper edge-coloring. Also,   $\chi'_{\rm {inj}}(G)$
may significantly differ from the injective chromatic number of the line graph $L(G)$ of $G$. In fact, if   $G$ has no vertices of degree $2$,  the 
 injective chromatic number of  $L(G)$ equals the strong chromatic index of $G$.
Recall that the {\em strong chromatic index}, {$\chi'_s(G)$},
is the minimum $k$ such that one can color the edges of $G$ with $k$ colors so that every two edges at distance at most $1$ have distinct colors.
By definition, $\chi_{\rm inj}'(G)\leq {\chi'_s(G)}$ for every graph $G$ and the difference between them can be large.  For example,
for the star $K_{1,n}$ we have $\chi'_{\rm {inj}}(K_{1,n})=1$ and $\chi_{\rm {inj}}(L(K_{1,n}))={\chi'_s(K_{1,n})} =n$.

\subsection{Previous results}

Injective vertex coloring was introduced and studied by  Hahn, Kratochv\'{\i}l,    Sotteau and  \v{S}ir\'a\v{n}~\cite{hss} in 2002. In particular, 
they showed that for each graph $G$ with maximum degree $\Delta$,  
$$\Delta\leq \chi_{\rm {inj}}(G)\leq \Delta(\Delta-1)+1,$$ and both bounds are sharp.

The notion of injective edge coloring was introduced in 2015 
by Cardoso, Cerdeira,  Cruz, and  Dominic~\cite{CCCD} motivated by a Packet Radio Network problem and independently in 2019 by
Axenovich,  D\"orr, Rollin, and  Ueckerdt~\cite{axia2019} (they called it {\em  induced  star arboricity}).

Cardoso et al.~\cite{CCCD} proved that computing $\chi'_{\rm {inj}}(G)$ of a graph $G$ is NP-hard and determined the injective chromatic index for
paths, cycles, wheels,  Petersen graph and complete bipartite graphs. They also proved that $\chi'_{\rm {inj}}(T)\leq 3$ for each tree $T$ and that {$\chi_{\rm inj}'(K_{\Delta+1})=\frac{\Delta(\Delta+1)}{2}$}.
Axenovich et al.~\cite{axia2019} concentrated more on another parameter,  {\em  induced   arboricity}, but they also proved that the  induced  star
arboricity of each planar graph is at most $30$ and can be as large as $18$. Apart from this, they presented bounds on 
 the  induced  star
arboricity of a graph in terms of its acyclic chromatic number and treewidth.

Bu and Qi~\cite{BQ} gave upper bounds on injective chromatic index of graphs with maximum degree $3$ and $4$ and low maximum average degree. In particular,
they showed that the injective chromatic index of every subcubic graph with maximum average degree at most $\frac{18}{7}$ (respectively,
at most $\frac{5}{2}$) is at most $6$ (respectively,
at most ${5}$).

Ferdjallah,  Kerdjoudj and  Raspaud~\cite{FKR} { used  Proposition 2.2 of \cite{CCCD}} to observe that the  induced  star
arboricity of a graph is exactly its injective chromatic index and proved a series of bounds on injective chromatic index of ``sparse" graphs.
They proved that for every $\Delta\geq 3$ and any graph $G$ with maximum degree at most $\Delta$,
\begin{equation}\label{D1}
{\chi_{\rm inj}'(G)\leq 2(\Delta- 1)^2.}
 \end{equation}


Ferdjallah,  Kerdjoudj and  Raspaud~\cite{FKR} posed the following conjecture.

\begin{conjecture}\label{CO1}
For every subcubic graph $G$,
$\chi_{\rm inj}'(G)\leq 6.$
\end{conjecture}

Furthermore, for bipartite graphs, Ferdjallah et al.~\cite{FKR} proved stronger bounds: for any bipartite graph $G$ with maximum degree at most $\Delta$,
\begin{equation}\label{D2}
 {\chi_{\rm inj}'(G)\leq \Delta(\Delta- 1)+1,}
 \end{equation}
and for every subcubic bipartite graph $G$,  $\chi_{\rm inj}'(G)\leq 6.$ 
They posed the following conjecture: \begin{conjecture}\label{CJBIP}
For every subcubic bipartite  graph, $\chi_{\rm inj}'(G)\leq 5$. 
\end{conjecture}
 If this conjecture is true, then $5$ is tight (see\cite{FKR}).
 
Ferdjallah et al.~\cite{FKR} also gave the exact upper bound of $5$ for the
 injective chromatic index of subcubic outerplanar graphs, and somewhat strengthened the result of Bu and Qi~\cite{BQ} for subcubic graphs mentioned above.

For graphs without $4$-cycles Mahdian~\cite{Ma00} proved  stronger upper bounds even for strong chromatic index:

\begin{theorem}[Mahdian~\cite{Ma00}]\label{Mah}
For every $C_4$-free graph $G$, $$\quad\chi_{s}'(G) \leq (2+o(1))\frac{\Delta^2}{\ln \Delta}.$$
\end{theorem}

The strong chromatic index of bipartite graphs was also studied. In particular,
 Faudree, Gy\'arf\'as, Schelp   and  Tuza~\cite{FGST90} in 1990 conjectured that $\chi_{s}'(G) \leq \Delta^2$
for every bipartite graph $G$ with maximum degree $\Delta$.
  Brualdi and Massey~\cite{BM93} posed in 1993 the refined conjecture that $\chi_{s}'(G) \leq \Delta(X)\Delta(Y)$
 for  every bipartite graph $G$ with parts $X$ and $Y$, where $\Delta(X)$ (resp. $\Delta(Y)$) is the maximum degree of a vertex of $X$ (resp. $Y$). Partial cases of this conjecture were proved by
Brualdi and Massey \cite{BM93},
  Nakprasit \cite{NAK2008} and Huang,  Yu and  Zhou~\cite{HYZ2017}.

\subsection{Our results}
The goal of this paper is to present new bounds on the  injective chromatic index of graphs with given maximum degree involving chromatic number.

Our first results are two steps toward Conjecture~\ref{CO1}:

\begin{theorem}\label{cubic}
For every subcubic graph $G$, $\quad\chi_{\rm {inj}}'(G) \leq 7$.
\end{theorem}

\begin{theorem}\label{planar}
For every planar subcubic graph $G$, $\quad\chi_{\rm {inj}}'(G) \leq 6$.
\end{theorem}

The bound in Theorem~\ref{planar} is exact: it is attained at $K_4$ and the $3$-prism.

The proof of this theorem yields a stronger bound for bipartite graphs which is a step forward Conjecture \ref{CJBIP}: it
implies that $\;\chi_{\rm {inj}}'(G) \leq 4$ for every bipartite planar subcubic graph. This bound is attained at the $3$-dimensional cube $Q_3$ with any edge deleted.

The main result of this paper is the following bound significantly improving~\eqref{D2}:

\begin{theorem}\label{main}
Let $2\leq \chi\leq \Delta$. If $G$ is a graph with maximum degree $\Delta$ and chromatic number $\chi$, then
\begin{equation}\label{me1}
\chi_{\rm {inj}}'(G) \leq(\chi-1)  \lceil 27 \Delta \ln \Delta\rceil.
\end{equation}
In particular, if $G$ is bipartite, then $\chi_{\rm {inj}}'(G) \leq \lceil 27 \Delta \ln \Delta\rceil.$
\end{theorem}

We also discuss the bound of Theorem~\ref{main} and compare it with similar bounds for $\chi'_s$.
First, we show that without restrictions on the chromatic number, the bound is much weaker even for graphs with large girth.
In fact, the order of magnitude of the bound in Theorem~\ref{Mah} cannot be improved not only for strong chromatic index but also for injective chromatic index even for graphs with arbitrary girth:

\begin{prop}\label{pr1}
For  every  $\Delta\geq 3$ and $g\geq 3$, there exists a   graph $G$ with maximum degree $\Delta$ and girth at least $g$ such that $\chi_{\rm {inj}}'(G) \geq\frac{\Delta(\Delta-1)}{4\ln \Delta }$.
\end{prop}


Then we show that the bound in Theorem~\ref{main} cannot be made less than $\Delta$ even for bipartite graphs with any girth.

\begin{prop}\label{pr2}
{For every $\Delta$-regular bipartie graph $G$,  $\chi_{\rm {inj}}'(G) \geq {\Delta }$.}
\end{prop}

Finally, we show that the result of Theorem~\ref{main} does not hold for strong chromatic index. Moreover, we show that there are bipartite 
$\Delta$-regular graphs of large girth that do not have ``too large" induced matchings.

\begin{prop}\label{pr3}
For  every  $\Delta\geq 21$ and $g\geq 3$, there exists a $\Delta$-regular bipartite graph $G$ with  girth at least $g$ such that 
the size of every induced matching in $G$ is less than 
$k=k(\Delta,n)=\left\lceil \frac{3n\ln \Delta}{\Delta}\right\rceil,$
where $n$ is the number of vertices in each of the parts of $G$. In particular, for each such $G$,
{$$\chi_{s}'(G) \geq  \frac{\Delta^2}{3\ln \Delta}.$$}
\end{prop}


In the next section, we prove Theorems~\ref{cubic} and~\ref{planar},  in  Section~3   prove the main result, Theorem~\ref{main}, and
in Section~4 
  Propositions~\ref{pr1}--\ref{pr3} are proven.

\section{Subcubic graphs}

\subsection{A bound for all subcubic graphs}
Another way to define the injective chromatic index of a graph $G$ is to consider the graph $G^{(\ast)}$ obtained from $G$ as follows: $V(G^{(\ast)})=E(G)$ and two vertices of $G^{(\ast)}$ are adjacent if the edges of $G$ corresponding to these two vertices of $G^{(\ast)}$ are at distance {1} in $G$ or in a triangle. 



Then
\begin{equation}\label{inj-cubic}
\chi_{\rm {inj}}'(G)=\chi(G^{(\ast)}).
\end{equation}

 Also, 
\begin{equation}\label{deg**}
\mbox{$\Delta({G^{(\ast)}})\leq 2(\Delta-1)^2$.}
\end{equation}

We will apply to $G^{(\ast)}$ the following theorem of Lov\'asz~\cite{Lo66}.

\begin{theorem}[Lov\'asz \cite{Lo66}]\label{tLov}
Let $G$ be a multigraph with maximum degree $\Delta$.
Let $t,k_1,k_2, \cdots, k_t$ be nonnegative integers such that $$ k_1+k_2+ \cdots + k_t=\Delta  -t +1.$$
Then the vertices of $G$ can be partitioned into sets $V_1, V_2, \cdots V_t$ so that the subgraph induced by each $V_i$ has maximum degree at most $k_i$.
\end{theorem}

Our goal is to prove Theorem~\ref{cubic}: {\em For every subcubic graph $G$, $\chi_{\rm {inj}}'(G) \leq 7$.}
We will study minimum counterexamples to the theorem and will show that they do not exist. We will exploit {\em partial injective $7$-colorings of edges} of graphs when not every edge is colored. Given a partial injective edge coloring $f$ of $H$ with colors from $[7]$ and an uncolored edge $e\in E(H)$,
$C_f(e)$ denotes the set of colors in $[7]$ not used on the edges at distance $1$ from $e$ or in a common triangle with $e$.
Furthermore,  $\overline{C}_f(e)$ stands for $[7]\setminus C_f(e)$.

We will use only Claims (d) and (f) of the   lemma  below, but it will be convenient to prove all of them in alphabetical order.

\begin{lemma}\label{co1}
Let $H$  be a counterexample to Theorem~\ref{cubic} with minimum $|E(H)| + |V (H)|$. Then\\
(a) $H$ is connected and has at least $8$ edges;\\
(b) $H$ is $3$-regular;\\
(c) $H$ does not contain a $K_4-e$;\\
(d) $H$ does not contain a triangle;\\
(e) $H$ does not contain a $K_{2,3}$;\\
(f) $H$ does not contain a $4$-cycle.
\end{lemma}

{\bf Proof.} Claim (a) immediately follows from the minimality of $H$ and the fact that we have $7$ available colors.

\smallskip
Suppose $H$ has a vertex $v$ of degree at most $2$. The case when $d(v)\leq 1$ is trivial, so suppose $N(v)=\{u,w\}$.
 By the minimality of $H$, graph $H'=H-v$ has an
injective edge coloring $f$ with colors in $[7]$. Since the degrees of $vu$ and $vw$ in $H^{(\ast)}$ are at most $6$, we can
greedily
 extend $f$ to  $vu$ and $vw$.  This contradicts the choice of $H$ and hence proves (b).

\smallskip
Suppose $H$ contains a copy $F$ of $K_4-e$ with $V(F)=\{a,b,c,d\}$. If {$H[V(F)]=K_4$}, then either $H$ is disconnected or $H=K_4$ and has $6$ edges.
Both possibilities contradict (a). Thus we may assume $bd\notin E(H)$. By the minimality of $H$, $H' = H-\{a,c\}$ has an injective edge coloring $f$ using colors in $[7]$. Each $e\in \{ab, bc, cd, da\}$ has at most three colored edges in $H'$ at distance one, so $|C_f(e)|\geq 4$. Furthermore, $|C_f(ac)|\geq 5$.
Thus we can color greedily these edges in the order $ab,bc,cd,da,ac$. This proves (c).

\smallskip
Suppose  $H$ contains a copy $F$ of $K_3$ with $V(F)=\{u,v,w\}$.
Let $u', v',$ and $w'$ denote the neighbor of $u,v$ and $w$ outside $F$, respectively. By (c), $u', v'$ and $w'$ are pairwise distinct. 
By the minimality of $H$,  $H' = H \setminus V(F)$ has an injective edge coloring $f$ using colors in $[7]$.
For each $e$ in $E_0=\{uv,vw,wu,uu',vv',ww'\}$, we have $|C_f(e)|\geq 3$. The maximum degree of $H^{(\ast)}[E_0]$ is $3$, and
$H^{(\ast)}[E_0]$ does not contain $K_4$. So by the list version of Brooks' Theorem (see e.g.~\cite{Viz} or~\cite{ERT})
$H^{(\ast)}[E_0]$ 
 is $3$-choosable. Thus we can extend  $f$ to whole $H$. This contradiction proves (d).

\smallskip
Suppose  $H$ contains a copy $F$ of $K_{2,3}$ with parts $\{x,y,z\}$ and  $\{u,v\}$. By (d), $H[V(F)]=F$.
Let $x',y'$ and $z'$ be the neighbor of $x,y$ and $z$, respectively, not in $F$. Some of them may coincide.
By the minimality of $H$,  $H' = H \setminus \{u,v\}$ has an injective edge coloring $f$ using colors in $[7]$.
For each $e$ in $E_0=\{xu,xv,yu,yv,zu,zv\}$, we have $|C_f(e)|\geq 3$: for example, the colored neighbors of $xu$ are $yy',zz'$ and 
two edges incident to $x'$. 
{
The maximum degree of $H^{(\ast)}[E_0]$ is $2$. {So by the list version of Brooks' Theorem, graph 
$H^{(\ast)}[E_0]$ 
 is $3$-choosable}. Thus we can extend  $f$ to whole $H$. This contradiction proves (e).}

\smallskip
Finally, suppose  $H$ has a $4$-cycle {$C=wxyzw$}. Let $w',x',y'$ and $z'$ be the neighbor of $w,x,y$ and $z$ outside of $C$.
 By (d), $C$ has no chords. By (e), all $w',x',y',z'$ are distinct.
 By the minimality of $H$,  $H' = H \setminus \{w,x,y,z\}$ has an injective edge coloring $f$ using colors in $[7]$.
 Let  $E_0=\{ww', xx', yy', zz'\}$  and $E_1=\{wx,xy,yz,zw\}$. Since each $e\in E_0$ has
  at most four colored edges at distance one, $|C_f(e)|\geq 3$.

\smallskip
{\bf Case 1:} $C_f(ww')\cap C_f(yy')\neq \varnothing$, say
 $1\in C_f(ww')\cap C_f(yy')$. Extend $f$ to $ww', yy'$ by $f(ww') = f(yy') = 1$. Choose distinct $f(xx')\in C_f(xx')$ and  $f(zz')\in C_f(zz')$ from the colors available for them. By symmetry, assume $f(xx') = 2, f(zz') = 3$. 
 
 For each $e\in E_1$,  $C_f(e)\neq \emptyset$ and $e$ has only one neighbor in $H^{(\ast)}[E_1]$. So, if we cannot 
 extend to $E_1$, then by symmetry we may assume 
 \begin{equation}\label{au1}
 C_f(wx) = C_f(zy) = \{7\}.
 \end{equation}
  Denote the set of edges incident to $w'$ except $ww'$ by $U(w')$, and define $U(x'), U(y'), U(z')$ similarly. Then in order for~\eqref{au1} to hold, we need 
  $f(U(w')\cup U(x')) = \{2,4,5,6\}$ and $f(U(z')\cup U(y')) = \{3,4,5,6\}$. Now since $f(ww') = f(yy')$, there are at least two available colors for $xx'$. Let $\alpha$ be such color for $xx'$ distinct from $2$.  Recolor $xx'$ by $\alpha$. Now both $2$ and $7$ are available for $zy$. Let $f(zy) = 2, f(wx) = 7$. 

If $2\in f(U(w'))$, then $2\notin f(U(x'))$. We can then extend $f$ by letting $f(xy) = 2, f(wz) = 7$, if $\alpha \neq 7$. 
If  $\alpha = 7$, let  $f(xy) = 7$.  Now for $wz$ we have at most 6 forbidden  colors, so we can extend $f$ to $wz$, as well.

If $2\in f(U(x'))$ then $2\notin f(U(w'))$. In this case, we extend $f$ letting $f(xy) = 7$ and $f(wz) = 2$.
In either case $f$ is an injective edge coloring of $H$, a contradiction.

\smallskip
{\bf Case 2:}
 $C_f(ww')\cap C_f(yy')= \varnothing =C_f(xx')\cap C_f(zz')$. Since $|C_f(e)|\geq 3$ for each $e\in E_0$, we can extend $f$ to the edges in $E_0$ so that
 $f(xx')\neq f(yy')\neq f(zz')\neq f(ww')\neq f(xx')$. By the case, the colors of all edges in $E_0$ are distinct.
  By symmetry, we may assume $f(ww') = 1, f(xx') = 2, f(yy') = 3, f(zz') = 4$, and $C_f(wx) = C_f(zy) = \{7\}$. Then $f(U(w')\cup U(x')) = \{1,2,5,6\}$ and
   $f(U(z')\cup U(y')) = \{3,4,5,6\}$. If we could recolor $xx'$ with some $\alpha\in C_f(xx')- \{1,2,3\}$, 
  then we  let $f(xx') = \alpha, f(zy) = 2, f(wx) = 7$. As in Case 1, depending on whether $2\in f(U(w'))$ or $2\in f(U(x'))$, we can either color $wz$ by $7$ and $xy$ by $2$ or vice versa. 
 
 Hence there is no such $\alpha$. Again by symmetry we may assume that
  initially  $C_f(ww') = \{1,2,4\}, C_f(xx') = \{1,2,3\}, C_f(yy') = \{2,3,4\}$, and $C_f(zz') = \{1,3,4\}$. But this is a contradiction to $C_f(ww')\cap C_f(yy') = \varnothing$.
\qed

 {\bf Proof of Theorem~\ref{cubic}.} 
{Let $H$ be a minimal counterexample as in Lemma \ref{co1}.}
 By Theorem~\ref{tLov},  the set of vertices of $H^{(\ast)}$ can be partitioned into two sets $V_1$ and $V_2$ so that  $\Delta(H^{(\ast)}[V_1])\leq 3$ and 
 $\Delta(H^{(\ast)}[V_2])\leq 4$.
 If $\chi(H^{(\ast)}[V_1])\leq 3$ and $\chi(H^{(\ast)}[V_2])\leq 4$, then
 we are done. By Brooks' Theorem, if $\chi(H^{(\ast)}[V_1])\geq 4$, then $H^{(\ast)}[V_1]$ contains a $K_4$,
  and if $\chi(H^{(\ast)}[V_2])\geq 5$, then $H^{(\ast)}[V_2]$ contains a $K_5$. So, we have two cases.
  
  {\bf Case 1:} $H^{(\ast)}[V_1]$ contains a $K_4$. 
  Let $e_1,e_2,e_3,e_4$ be the vertices of this $K_4$ and let $e_i=v_{2i-1}v_{2i}$ for $i\in [4]$. 
  Since $H$ has no $3$-cycles, all $v_1,\ldots,v_8$ are distinct and all edges $e_1,e_2,e_3,e_4$ are at distance exactly $1$ from each other in $H$.
  By symmetry, we may assume that $v_1$ is adjacent   to $v_3$ and $v_5$, and $v_2$ is adjacent to $v_7$.
  Then since $H$ has no $3$- and $4$-cycles, in order to have $e_2$ and $e_3$ at distance $1$, we need $v_4v_6\in E(H)$.
  Now for the same reason, $v_7$ is not adjacent to $v_3$ or $v_5$ and
  neither of $v_7$ and $v_8$ can be adjacent to two vertices in the $5$-cycle $v_1v_3v_4v_6v_5v_1$. So again by symmetry, we may assume 
  that $v_7$ is adjacent to $v_4$, and $v_8$ is adjacent to $v_5$ (see the picture below).
    
 \begin{center}
\includegraphics[scale=0.4]{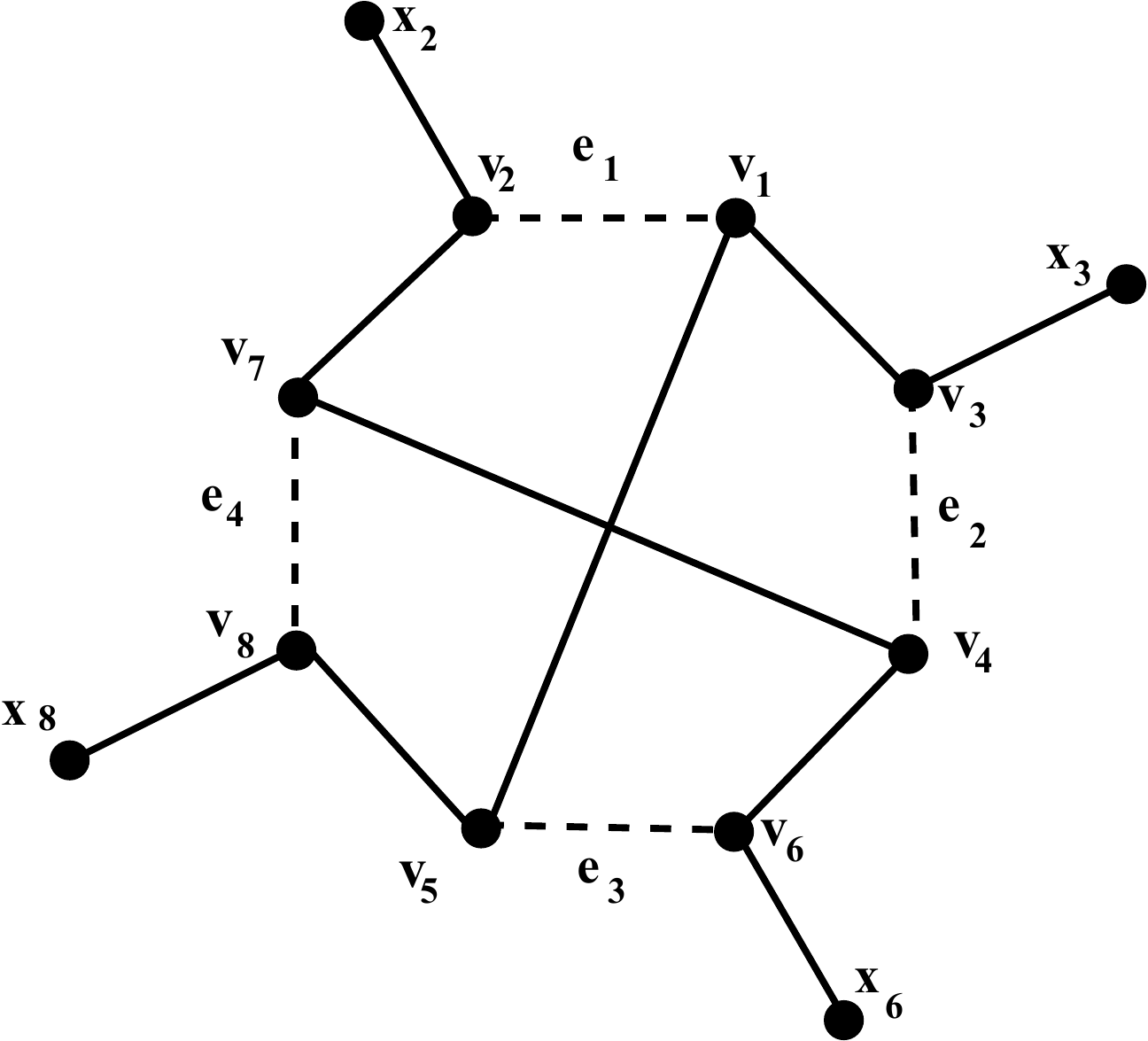}
\end{center}
Since $\Delta(H^{(\ast)}[V_1])\leq 3$, all the edges incident to $\{v_1,\ldots,v_8\}$ apart from $e_1,e_2,e_3$ and $e_4$  are vertices in $V_2$ in
$H^{(\ast)}$. In particular, vertex $v_1v_5\in V_2$ is adjacent in $H^{(\ast)}$ to  vertices $v_2x_2,v_2v_7,v_3x_3,v_8x_8,$ $v_6x_6,v_6v_4$ in $V_2$,
contradicting   $\Delta(H^{(\ast)}[V_2])\leq 4$.

{\bf Case 2:} $H^{(\ast)}[V_2]$ contains a $K_5$.   
Let $e_1,\ldots,e_5$ be the vertices of this $K_5$ and let $e_i=v_{2i-1}v_{2i}$ for $i\in [5]$. Let $V_0=\{v_1,\ldots,v_{10}\}$.
As in Case 1, all $v_1,\ldots,v_{10}$ are distinct and all edges $e_1,\ldots,e_5$ are at distance exactly $1$ from each other in $H$.
In order to achieve this, each vertex in $V_0$ has neighbors only in $V_0$. So by the minimality of $H$, $|V(H)|=10$. The only cubic 
$10$-vertex graph with  no $3$- and $4$-cycles is
  the Petersen graph ${\bf P}$ {(see \cite{BIG1993} or  \cite{HS1993})}. As it was mentioned in the introduction,  Cardoso et al.~\cite{CCCD} proved that
   $\chi_{\rm {inj}}'({\bf P})=5$. Hence $H^{(\ast)}[V_2]$ is not $K_5$, a contradiction again.
\qed
 

\subsection{A bound for planar subcubic graphs}

We will use the following partial case of a generalization of Brooks' Theorem proved independently by Borodin~\cite{Bor} and 
Bollob\' as and Manvel~\cite{BMa}.

\begin{theorem}[\cite{BMa,Bor}]\label{BorBol}
If $G$ is a connected graph with $\Delta(G)\leq 3$ not containing $K_4$, 
then one can partition $V(G)$ into sets $X$ and $Y$ so that $X$ is independent and $G[Y]$ is an acyclic graph with maximum degree at most $2$.
\end{theorem}

 {\bf Proof of Theorem~\ref{planar}.} Let $G$ be a vertex-minimal subcubic plane graph with $\chi'_{\rm {inj}}(G)\geq 7$.
 Then, $\delta(G)\geq 2$, $G$ is connected,  and has at least $7$ edges; in particular, $G\neq K_4$. So
 by Theorem~\ref{BorBol}, $V(G)=X\cup Y$ where $X$ is independent and $G[Y]$ is an acyclic graph with maximum degree at most $2$.

Construct an auxiliary (multi)graph $G'$ with vertex set $Y$ as follows. For each $x\in X$, if $N(x)=\{y_1,y_2\}$, then delete $x$ and add edge $y_1y_2$,
and if $N(x)=\{y_1,y_2,y_3\}$, then delete $x$ and add edges $y_1y_2$, $y_1y_3$
and $y_2y_3$.  By construction,  $G'$ is a planar (multi)graph. 
By the Four Color Theorem, $G'$ has a 4-coloring $g$.

We now color $E(G)$ in two steps.  On Step 1   for each edge $xy$ connecting $X$ with $Y$ color $xy$ with $g(y)$.
By construction, for each $x\in X$, the colors of all edges incident with $x$ are distinct. Also, two edges of the same color cannot have an edge in $Y$ connecting them. So, after Step 1, we  have a partial injective edge coloring of $G$  with colors in $[4]$.

On Step 2  for each path component $y_1,y_2,….,y_s$ in $G[Y]$, color the first two edges with $5$, second pair of edges with $6$, third pair again with $5$ and so on. This yields an injective edge-coloring of $G$.\qed

\medskip
If our planar subcubic graph $G$ is also bipartite, then instead of the partition $V(G)=X\cup Y$ provided by Theorem~\ref{BorBol}, we can use the natural bipartition of $G$, and do not need to run Step 2 and use the extra colors $5$ and $6$. Thus our proof has the following implication.

\begin{corollary}\label{Andre}
For every planar subcubic bipartite graph $G$, $\quad\chi_{\rm {inj}}'(G) \leq 4$.
\end{corollary}

The bound in the corollary is sharp: $\chi_{\rm {inj}}'(Q_3)=4$ where $Q_3$ is the graph of the unit $3$-cube. Also, if we delete any edge from $Q_3$,
the injective chromatic index of the remaining graph is still $4$.

\section{Graphs with high maximum degree}

Our  tool in this section is  Lov\' asz Local Lemma~\cite{EL} in a slightly stronger form proved by Spencer~\cite{epd}:

\begin{theorem}[Lov\' asz Local Lemma~\cite{EL,epd}]\label{LLL}
Let $A_1, \dots,A_n$ be events such that $Pr[A_i]\leq p$, for $1\leq i\leq n$. Suppose each event is independent of all the other events except for at most $d$ of them. If $ep(d+1)< 1$, then $Pr[\bigwedge_{i = 1}^n \overline{A_i}]>0$.
\end{theorem}

A subset $F$ of edges of a graph $G$ is $G$-{\em good} if no two edges in $F$ are at distance $1$ in $G$ or in the same triangle.
In other words, a $G$-{good} set is an independent set in the auxiliary graph $G^{(\ast)}$ defined in Section~3, so that an injective edge coloring is simply a partition of $E(G)$ into $G$-{good} sets.

For $\Delta\leq 40$ the bound of Theorem~\ref{main} is weaker than~\eqref{D1}. So, it is enough to prove Theorem~\ref{main} for $\Delta\geq 41$. If $\chi(G)=\chi$, then there is a coloring of $G$ with color classes $Y_1,Y_2,\ldots,Y_\chi$ such that for every
$1\leq j\leq \chi-1$, set $Y_j$ is a maximal (by inclusion) independent set in the graph $G[Y_j\cup Y_{j+1}\cup\ldots \cup Y_\chi]$.
The theorem below allows us for each $1\leq j\leq \chi-1$ to color the edges connecting $Y_j$ to 
$G[Y_{j+1}\cup Y_{j+2}\cup\ldots \cup Y_\chi]$ with $\lceil 27  \Delta \ln \Delta \rceil$ colors.
Thus, applying this theorem $\chi -1$ times will imply  Theorem~\ref{main}.

  
\begin{theorem}\label{main2}
Let $ \Delta\geq 41$ and $k=\lceil 27  \Delta \ln \Delta \rceil$. If $G$ is a graph with maximum degree $\Delta$ and $Y$ is a {maximal} independent set in $G$,
then we can partition $E(Y,V(G)-Y)$ into
  $k$ $G$-{good} sets.
\end{theorem}

{\bf Proof of Theorem~\ref{main2}.}
Let $X=N(Y)=V(G)-Y$. We will construct $k$ random $G$-{good} sets $J_1,\ldots,J_k$,  using the following algorithm for $j=1,\ldots,k$.

\medskip
{\bf Step 1:} Construct a random subset $X_j$ of $X$ by including each vertex of $X$ into $X_j$ with probability $\frac{1}{\Delta}$ independently of each other.

{\bf Step 2:} Delete from  $ X_j$ every vertex that after Step 1 has a neighbor in  $X_j$. 

{\bf Step 3:} Let $J_j$ be the set of the edges connecting $X_j$ with $Y$.

{\bf Step 4:} For each $y\in Y$, if $y$ had at least two neighbors in $X_j$ {\em after Step 1}, then
 remove all these edges from $J_j$ incident to $y$. 

\medskip
We claim that 
\begin{equation}\label{e3}
\mbox{$J_j$ is $G$-good.}
\end{equation}
 Indeed, since after Step 2, $X_j$ is independent, no two edges in $J_j$ are in a common triangle. Suppose that  edges $x_1y_1$ and 
$x_2y_2$ in $J_j$ are at distance $1$, say  $x_1y_2\in E(G)$. But then after Step $1$ vertex $y_2$ has at least  two neighbors in $X_j$ and hence on Step 4 edges $y_2x_1$ and $y_2x_2$  would be deleted from $J_j$. This proves~\eqref{e3}.

By~\eqref{e3}, it remains to prove that with positive probability each edge in $E(X,Y)$ will belong to at least one $J_j$.
For this, we introduce several events and estimate their probabilities. 

\medskip
 Denote the event that an {\em $x\in X$ is in $X_j$ after Step 1} by $F_{1,j}(x)$.  By definition, $p[F_{1,j}(x)] = 1/\Delta$ for each $x\in X$.

Denote the event that an {\em $x\in X$ is in $X_j$ after Step 2}
 by $F_{2,j}(x)$.  Then
for each $x\in X$,
\begin{equation}\label{e1}
\frac{1}{\Delta}=p[F_{1,j}(x)]\geq 
p[F_{2,j}(x)] = \frac{1}{\Delta}\cdot \left( \frac{\Delta-1}{\Delta}\right)^{|N(x)\cap X|}\geq \frac{1}{\Delta}\cdot \left( \frac{\Delta-1}{\Delta}\right)^{\Delta-1}>\frac{1}{e\cdot \Delta}.
\end{equation}


For $xy\in E(X,Y)$ let  $F_j(xy)$ be the event that {\em $xy$ is in $J_j$ after Step 4}.

\medskip
Observe that for each $xy\in E(X,Y)$,
$$
p[F_j(xy)] = p[F_{2,j}(x)]\prod_{x'\in N(y)-x} (1-p[F_{1,j}(x')]).$$
Hence by~\eqref{e1},
\begin{equation}\label{e2}
p[F_j(xy)] \geq 
 \frac{1}{e \Delta}\left(1-\frac{1}{\Delta} \right)^{d(y)-1}\geq \frac{1}{e \Delta}\left(1-\frac{1}{\Delta} \right)^{\Delta-1}
  >\frac{1}{e^2\Delta}.
\end{equation}

For $xy\in E(X,Y)$, let $A(xy)$ denote the event that none of $F_j(xy)$ happened.
We want to show that with positive probability none of $A(xy)$ occurs, because in this case we can assign to each edge one of the $k$
colors. We plan to apply Theorem~\ref{LLL}, so we need to give upper bounds on the probability of each $A(xy)$ and on the number of events $A(x'y')$ on which depends $A(xy)$.

The first part is easy,
 since for distinct $j$ the events $F_j(xy)$ are independent, and hence by~\eqref{e2}
\begin{equation}\label{e4}
   p[A(xy)] =\prod_{j=1}^k (1-p[F_j(xy)])\leq 
   \left(1-\frac{1}{e^2\Delta}\right) ^{27  \Delta \ln \Delta} <\exp\{{-\frac{27 \ln \Delta}{e^2}}\}=\Delta^{-{27}/{e^2}}.
   \end{equation}

For the second part, we define $Q(xy)$ as the set of edges $x'y'$ such that the distance from $\{x',y'\}$ to
$N(y)\cup (N(x)\cap X)$ is at most $1$ and will prove that for each $xy\in E(X,Y)$
\begin{equation}\label{28}
\mbox{\em $A(xy)$ is independent of all $A(x'y')$ such that  $x'y'\in E(X,Y)-Q(xy)$.}
\end{equation}

Recall that $A(xy)$ is fully defined by the events $F_1(xy),F_2(xy),\ldots,F_k(xy)$, and that for all $j'\neq j$
event $F_j(xy)$ is independent of all $F_{j'}(x'y')$. Thus to get~\eqref{28}, it is enough to
prove that for each $xy\in E(X,Y)$ and each $j\in [k]$,
\begin{equation}\label{281}
\mbox{\em $F_j(xy)$ is independent of all $F_j(x'y')$ such that  $x'y'\in E(G)-Q(xy)$.}
\end{equation}

By definition, for each $xy\in E(X,Y)$, event $F_j(xy)$ occurs if and only if $F_{1,j}(x)$ occurs but none of 
${F_{1,j}(x')}$ over all $x'\in (N(x)\cap X)\cup (N(y)-x)$ occurs.
 Since the events   $F_{1,j}(x)$  are independent for all $x\in X$,  $F_j(xy)$ is fully defined by the event
$$R(xy)=\bigcup\nolimits_{x'\in N(y)\cup   (N(x)\cap X)   }F_{1,j}(x').$$

It follows that $F_j(xy)$ is independent of all $F_j(x_1,y_1)$ such that $R(xy)\cap R(x_1,y_1)=\emptyset$.
Since  each $F_{1,j}(x')$ belongs to $R(xy)$ only for the edges $xy$ such that
 $\{x,y\}$ is at distance at most $1$ from $x'$,~\eqref{281} follows. This in turn yields~\eqref{28}.
 
 Since for each $x'\in X$
 there are at most $\Delta^2$  edges at distance at most $1$ from $x'$,~\eqref{28} implies
 that $F_j(xy)$   is independent of all but at most
 $$\left|\bigcup_{x'\in N(y)}F_{1,j}(x')\cup \bigcup_{x''\in N(x)\cap X}F_{1,j}(x'')\right| \Delta^2\leq (2\Delta-1) \Delta^2$$ 
 other events  $F_j(x'y')$.
 
   Therefore, by Theorem~\ref{LLL} with $d< 2 \Delta^3$ and
 $p=\Delta^{-{27}/{e^2}}$, it is enough to prove that 
 \begin{equation}\label{e5}
   e\cdot  2 \Delta^3\cdot \Delta^{-{27}/{e^2}}<1.
   \end{equation} 
To derive~\eqref{e5}, it is enough to check that $\Delta^{-3+{27}/{e^2}}>2e$.   Since $-3+{27}/{e^2}>0.654$, this holds for $\Delta>40$.\qed

\section{Lower bounds}

\subsection{Proofs of Propositions~\ref{pr1} and~\ref{pr2}}
We will use the following result obtained by Kostochka and Mazurova~\cite{KM} and independently by Bollob\' as~\cite{Bol81}:

\begin{theorem}[\cite{KM,Bol81}]\label{Bol}
For every  $\Delta $ and $g$, there is a graph $G_{\Delta,g}$ with maximum degree $\Delta$ of girth at least $g$ such that
$\alpha(G_{\Delta,g})\leq \frac{2\ln \Delta}{\Delta}|V(G)|$ and {$|E(G_{\Delta,g})|>(\Delta-1)|V(G)|/2$}.
\end{theorem}
  
 
\medskip
{\bf Proof of Proposition~\ref{pr1}.} Fix any $\Delta\geq 3$ and $g\geq 3$. Let $G_{\Delta,g}$ be a graph satisfying  Theorem~\ref{Bol}. Let $n=|V(G_{\Delta,g})|$.
Assume that $G$ has  an injective edge coloring with $k$ colors. Let $\{I_1,\cdots,I_k\}$ be a partition of $E(G)$ into  $k$ color classes, and let
$I_1$ be a largest color class. Then 
\begin{equation}\label{lower}
|I_1|\geq \frac{(\Delta-1) n}{2k}.
\end{equation} 

Let $V(I_1)$ denote the union of the vertex sets of all edges in $I_1$. By definition, every component of  $G[V(I_1)]$ is a star. Hence, if we delete 
from $V(I_1)$
the center of each star (when the star has only two vertices, then we assign exactly one of them as the center), then we obtain an independent set 
$J(I_1)$ of vertices in $G_{\Delta,g}$ with $|J(I_1)|=|I_1|$. Thus by~\eqref{lower} and the choice of $G_{\Delta,g}$,
$$ \frac{(\Delta-1) n}{2k}\leq \frac{2\ln \Delta}{\Delta}n.
$$
This yields the proposition.\qed

\medskip
{\bf Proof of Proposition~\ref{pr2}.} 
The argument is similar to the proof above (and simpler). 
{For each $\Delta$, 
 let
$B$ be any $\Delta$-regular bipartite graph. 
 By Marriage Theorem,
$\alpha(B)=0.5 |V(B)|$. As in the proof of Proposition~\ref{pr1}, the number of edges in any color class of an injective edge coloring 
of $B$ is at most
$$\alpha(B)=0.5 |V(B)|=\frac{|E(B)|}{\Delta},$$
so we need at least $\Delta$ colors for injective edge coloring of $B$.\qed
 }

\subsection{The  Bipartite Configuration Model}
In order to prove Proposition~\ref{pr3}, we will use a bipartite version of the configuration model.
This model  in different versions is due to Bender and Canfield~\cite{BC1} and Bollob\'as~\cite{B2}. 
The bipartite version is considered in several papers. We follow the convention and use the results described in survey~\cite{W2}(Section 3.2) by Wormald.

Let $n$ and $D$ be  positive integers, and
\begin{equation}\label{s71}
\mbox{\em 
 $V_n=V_n(D)=\{v_1,\ldots,v_{nD}\}$ and $W_n=W_n(D)=\{w_1,\ldots,w_{nD}\}$ be disjoint sets.} 
\end{equation}

A {\em configuration} (of order $n$ and degree $D$) is 
a perfect matching from $V_n$ to $W_n$ (each edge has one end in $V_n$ and one in $W_n$). 
Let $\mathcal{F}_D(n)$  denote the collection of all $(D n)!$ such matchings.

For every $F\in\mathcal{F}_D(n)$ we define the $D$-regular bipartite multigraph 
 $\pi(F)$ with parts ${X_n}=\{x_1,\ldots,x_n\}$ and ${Y_n}=\{y_1,\ldots,y_n\}$ as follows: For every $j\in [n]$ we glue the $D$ vertices\\
 $v_{D\cdot (j-1)+1}, v_{D\cdot (j-1)+2},\ldots,v_{D\cdot j}$ into a new vertex $x_j$ and  the $D$ vertices 
 $w_{D\cdot (j-1)+1}, w_{D\cdot (j-1)+2},\ldots,w_{D\cdot j}$ into a new vertex $y_j$.


\begin{defn}
Let $\mathcal{G}_{D,g}(n)$ be the set of all $D$-regular  {bipartite graphs}  with parts $X_n=\{x_1,\ldots,x_n\}$ and $Y_n=\{y_1,\ldots,y_n\}$
 and girth at least $g$, and 
$\mathcal{G}'_{D,g}(n)=\{F\in \mathcal{F}_D(n)\,:\, \pi(F)\in \mathcal{G}_{D,g}(n)\}$.
\end{defn}


\medskip
 Bollob\'as~\cite{B2} and Wormald~\cite{W1} proved that for each fixed $g$ and $D$, there is $\epsilon(g,D)>0$ such that
   $$\frac{|\mathcal{G}'_{D,g}(n)|}{|\mathcal{F}_D(n)|}>\epsilon(g,D).$$
As discussed in~\cite{W1} and~\cite{BMc}, the same phenomenon holds for bipartite configurations. So, we will use the following fact:

\begin{theorem}\label{BW}
For each fixed $D,g \geq 3$, if a property holds   for $\pi(F)$ for almost all configurations  $F\in \mathcal{F}_D(n)$, then
it also holds for   $\pi(G)$ for almost all  $G\in \mathcal{G}'_{D,g}(n)$.
\end{theorem}

\subsection{Proof of Proposition~\ref{pr3}}

\begin{lemma}\label{l18}
For each fixed $D\geq 21$ and almost all configurations  $F\in \mathcal{F}_D(n)$,
$\pi(F)$ does not have induced matchings of size 
\begin{equation}\label{s7}
k=k(D)=\left\lceil \frac{3n\ln D}{D}\right\rceil.
\end{equation}
\end{lemma}

{\bf Proof.} Given $V_n$ and $W_n$ as in~\eqref{s71}, the number of matchings of size $k$ between $V_n$ and $W_n$
corresponding to matchings between $X_n$ and $Y_n$ defined above is
${n\choose k}^2 \cdot D^{2k}\cdot k!:
$ There are ${n\choose k}^2$ ways to choose the $k$-element subsets of $X_n$ and $Y_n$ joined by a matching, then there are
$D^{2k}$ ways to choose the vertices in $V_n$ and $W_n$ to represent the chosen $2k$ vertices from $X_n$ and $Y_n$, and finally
there are $k!$ ways to match the $k$ chosen vertices of $V_n$ with the $k$ chosen vertices in $W_n$.

For each such matching $M$, the number of configurations  $F\in \mathcal{F}_D(n)$ in which $M$ is an induced matching 
in  {the multigraph} $\pi(F)$ is exactly
$$\left(\prod_{j=1}^{k(D-1)}(D(n-k)+1-j)\right)(D(n-k))!.
$$	
Hence the portion of $F\in \mathcal{F}_D(n)$ such that $\pi(F)$  has at least one induced matching of size $k$ is at most
$${n\choose k}^2 { D^{2k}\cdot k!}\left(\prod_{j=1}^{k(D-1)}(D(n-k)+1-j)\right)\frac{(D(n-k))!}{(Dn)!}
$$
$$\leq \left(\frac{n^k}{k!}\right)^2\left( D^{2k}\cdot k!\right) \left(\prod_{j=1}^{k(D-1)}\frac{D(n-k)+1-j}{Dn+1-j}\right) \frac{1}{(D(n-k))^k}$$
$$\leq \frac{(nD)^{2k}}{k!}\cdot \frac{1}{(D(n-k))^k} \left(\frac{D(n-k)}{Dn}\right)^{k(D-1)} 
= \frac{(nD)^{k}}{k!}\cdot \left(\frac{n-k}{n}\right)^{k(D-2)} $$
$$\leq \left[ \frac{  nD\cdot e}{k}\left(1-\frac{k}{n}\right)^{D-2}\right]^k\leq  \left[ \frac{  nD}{k} \cdot e^{1-(D-2)k/n}\right]^k.
$$ 

By~\eqref{s7}, $k\geq \frac{3n\ln D}{D}$, so the last expression in the brackets is at most
$$\frac{  nD^2}{3n \ln D} \exp\left\{1-\frac{3(D-2) \ln D}{D}\right\}<  \frac{  D^2}{ \ln D} \exp\left\{-\frac{3(D-2) \ln D}{D}\right\}.$$

Since  $D\geq 21$, $D-2> 0.9D$ and $\ln D>3$. So
$$ \frac{  D^2}{ \ln D} \exp\left\{-\frac{3(D-2) \ln D}{D}\right\}<  \frac{  D^2}{ 3} \exp\left\{-2.7 \ln D\right\}<\frac{1}{3D^{0.7}}.
$$	
	
	It follows that the portion of $F\in \mathcal{F}_D(n)$ such that $\pi(F)$ has at least one induced matching of size $k$ is at most $(3D^{0.7})^{-k}\to_{n\to\infty}0$.
\qed

\bigskip
Now we are ready to prove Proposition~\ref{pr3}.
By Lemma~\ref{l18} together with Theorem~\ref{BW}, for every $\Delta\geq 21$ and $g\geq 3$, there is a $\Delta$-regular bigraph $G$ with girth at least $g$ with the maximum size of an induced matching less than $\frac{3|V(G)|\ln \Delta}{2\Delta}$. Then $\chi'_s(G)>\frac{\Delta^2}{3\ln \Delta}$.

\section{Concluding remarks}

1. Our proof of Theorem~\ref{main2} does not work for injective list edge-coloring. We do not know how to prove the list analog of this theorem.

\medskip\noindent
2. On the other hand, several parts of the proof of Theorem~\ref{cubic} do work for list coloring. 

\medskip\noindent
3. Recall that an $L(h,k)$-coloring of a graph $H$ is a coloring $f$ of the vertices of $H$ with colors $1,2,\ldots$ such that for every adjacent {vertices} $x,y\in V(G)$,
$|f(x)-f(y)|\geq h$ and for each $u,v\in V(G)$ at distance exactly $2$, $|f(u)-f(v)|\geq k$. Such colorings arose from several applications and attracted some attention, see survey~\cite{Cal}. In these terms, 
if a graph $G$ is triangle-free (in particular, if $G$ is bipartite), then each injective edge-coloring of $G$ corresponds to an $L(0,1)$-coloring of $L(G)$ and
vice versa.

\bigskip
{\bf Acknowledgment.} We thank Kathie Cameron, Dieter Rautenbach and Nick Wormald for helpful discussions.
We also thank both referees for helpful comments.

\end{document}